\documentclass[reqno,11pt]{amsart}
\usepackage[latin1]{inputenc}
\usepackage[T1]{fontenc}
\usepackage{lmodern}
\usepackage{a4wide}
\usepackage{amsmath,amsfonts,amsthm,amssymb}
\usepackage{hyperref}
\usepackage{mathtools}
\usepackage{color}
\usepackage{dsfont}
\usepackage{pgf,tikz}

\usepackage{mathrsfs}
\usepackage{enumerate}
\usepackage[shortlabels]{enumitem}
\everymath{\displaystyle}
\usetikzlibrary{arrows}

\newcommand{\supp}{{\mathrm{supp}}}

\renewcommand\leq{\leqslant}
\renewcommand\geq{\geqslant}

\newcommand\be{\begin{equation}}
\newcommand\ee{\end{equation}}

\theoremstyle{plain}
\newtheorem{theorem}{\bf Theorem}[section]

\newtheorem{remark}{\bf Remark}[section]
\newtheorem{definition}{\bf Definition}[section]

\numberwithin{equation}{section}

\title[Poincar\'{e}-type inequalities in Musielak spaces]{Poincar\'{e}-type inequalities in Musielak spaces}

\author[]{Ahmed Youssfi}
\address{Ahmed Youssfi\\
Laboratory of Engineering, Systems and Applications (LISA)\\
National School of Applied Sciences\\ 
Sidi Mohamed Ben Abdellah University\\
My Abdellah Avenue, Road Imouzer, P.O. Box 72 F\`es-Principale, 30 000 Fez, Morocco}
\email{address:ahmed.youssfi@gmail.com ; ahmed.youssfi@usmba.ac.ma}

\author[]{Youssef Ahmida}
\address{Youssef Ahmida\\
	Laboratory of Engineering, Systems and Applications (LISA)\\
	National School of Applied Sciences\\ 
	Sidi Mohamed Ben Abdellah University\\
	My Abdellah Avenue, Road Imouzer, P.O. Box 72 F\`es-Principale, 30 000 Fez, Morocco}
\email{youssef.ahmida@usmba.ac.ma}

\begin{document}

\maketitle
\begin{abstract}
In this paper we investigate Poincar\'e-type integral  inequalities
in the functional Musielak structure. We extend the ones already well known in Sobolev, Orlicz and variable exponent Sobolev spaces. We introduce conditions on the Musielak functions under which they hold. The identification with null trace functions space is given.
\end{abstract}




\section{Introduction and main results}
In the last two decades, there has been an increasingly interest in studying  Musielak spaces, particularly for the analysis of nonlinear partial differential equations with non-standard growth conditions which come from modelling modern materials such as non Newtonian fluids, see for instance \cite{gwiazda2,raj-ru1} and the references therein.
\par In \cite{HH3,HH4,MJ} there is a basic background on the Musielak spaces $L_M(\Omega)$ and the Musielak-Sobolev spaces $W^{m}L_M(\Omega)$. An interesting missing feature is the Poincar\'e-type inequalities (in norm or in integral forms) in the closed subspace $W_0^{m}L_M(\Omega)$ defined  as the closure of the set $C^\infty_0(\Omega)$ of compactly supported functions in $\Omega$ with respect to the weak-$\ast$ topology $\sigma(\Pi L_M,\Pi E_{M^\ast})$ in the Musielak-Orlicz space $W^{m}L_{M}(\Omega)$.
However, proving the Poincar\'e integral inequality for functions in $C^\infty_0(\Omega)$ and then extending it by a density argument (as is often done for a constant exponent) is not an easy task since the passage to the limits is not allowed because of the lack in general of density of smooth functions in $W_0^{m}L_M(\Omega)$ at least in the modular sense (see Definition \ref{def:mod-conv}).
This is mainly due to the fact that the shift operator is not acting in general on Musielak spaces unless some regularity conditions on the Musielak function $M$ are satisfied see \cite{AYGS2017,AY2018}.
\par In this paper, we are interested in the problem of Poincar\'e-type integral inequality in the Musielak spaces. Such integral inequality yields obviously the  Poincar\'e norm inequality. Precisely, we give
sufficient conditions on the $\Phi$-function $M$ for the following Poincar\'e-type inequality
$$
\int_\Omega \sum_{|\alpha|< m} M(x,|D^\alpha u(x)|)dx\leq \int_\Omega \sum_{|\alpha|= m} M(x,c|D^\alpha u(x)|)dx
$$
to hold for every $u\in W_0^{m}L_M(\Omega)$ where $c>0$ is a constant. We also get the same inequality in the subspace $W_0^{m}E_M(\Omega)$ under minimal assumptions.
\subsection{Poincar\'{e}-type inequalities : State of the art}
Let $\Omega$ be a bounded open subset of $\mathbf{R}^N$, $N\geq1$ and let $1\leq p<\infty$. The usual Sobolev spaces are denoted $W^{1,p}(\Omega)$ while by $W^{1,p}_0(\Omega)$ we denote the norm closure of $C^\infty_0(\Omega)$ in $W^{1,p}(\Omega)$. The classical Poincar\'e integral inequality asserts that
\begin{equation}\label{poincare1}
\int_\Omega|u(x)|^{p}dx\leq C(\Omega,p)\int_\Omega|\nabla u(x)|^{p}dx
\end{equation}
for every $u\in W^{1,p}_0(\Omega)$ where $C(\Omega,p)$ is a constant depending on $\Omega$ and $p$. In fact, this inequality remains valid if $\Omega$ is only bounded in one direction. Recalling here that when $\Omega$ is regular (see for instance \cite[Theorem 4.14]{book_Necas}) we have
\begin{equation}\label{tr}
W^{1,p}_0(\Omega)=\big\{u\in W^{1,p}(\Omega) : tr(u)=0 \mbox{ on } \partial\Omega\big\}
\end{equation}
and hence $W^{1,p}_0(\Omega)=W^{1,1}_0(\Omega)\cap W^{1,p}(\Omega)$.
\par Gossez \cite[Lemma 5.7]{GJP1} proved the existence of two constants $c_m>0$ and $c_{m,\Omega}>0$ such the following Orlicz version of Poincar\'e integral inequality
\begin{equation}\label{orliczpoincare}
\int_\Omega \sum_{|\alpha|< m}\varphi(|D^\alpha u(x)|)dx\leq c_m\int_\Omega \sum_{|\alpha|= m}\varphi(c_{m,\Omega}|D^\alpha u(x)|)dx,
\end{equation}
holds for every $u\in W_0^{m}L_{\varphi}(\Omega)$. Here, $W^{m}_{0}L_{\varphi}(\Omega)$ is defined as the closure of the set $C^\infty_0(\Omega)$ of compactly supported functions in $\Omega$ with respect to the weak-$\ast$ topology $\sigma(\Pi L_\varphi,\Pi E_{\varphi^\ast})$ in the Orlicz spaces $W^{m}L_{\varphi}(\Omega)$, where $\varphi$ and $\varphi^\ast$ form a pair of complementary $N$-functions, cf. \cite{AF}. Since no extra condition is assumed on $\varphi$, inequality (\ref{orliczpoincare}) proved in $W^{m}_{0}L_{\varphi}(\Omega)$ covers not only (\ref{poincare1}) but it remains valid for a wide class of Orlicz functions. In contrast to Sobolev spaces $W^{1,p}_0(\Omega)$, the introduction of the Orlicz spaces $W^{m}_{0}L_{\varphi}(\Omega)$, defined by mean of the weak-$\ast$ topology $\sigma(\Pi L_\varphi,\Pi E_{\varphi^\ast})$, seems to be more convenient and very interesting in the theory of existence of PDEs in nonreflexive functional spaces, since firstly the weak topology is not equivalent in general to the strong one
and secondly coarser topology has more compact sets than the strong one.
\par Unfortunately, in the framework of variable exponent spaces the situation is more complicated and more regularities on the exponent are needed.
In fact in the Sobolev space $W^{1,p(\cdot)}_0(\Omega)$, defined as the norm
closure of $C^\infty_0(\Omega)$ functions in $W^{1,p(\cdot)}(\Omega)$, the Poincar\'e norm inequality was first proved in the pioneering paper \cite[Theorem 3.10]{KR} written about variable exponent Sobolev spaces provided that the exponent $p(\cdot)$ is continuous on $\overline{\Omega}$ and then by using the approche based on the boundedness of the maximal operator on $L^{p(\cdot)}(\Omega)$, the authors in \cite[Theorem 6.21]{bookCF} proved the  Poincar\'e norm inequality
\begin{equation}\label{poincarenorm}
\|u\|_{L^{p(\cdot)}(\Omega)}\leq c(N,p(\cdot),\Omega)\|\nabla u\|_{L^{p(\cdot)}(\Omega)}
\end{equation}
for every $u\in W_0^{1,p(\cdot)}(\Omega)$ and for exponents $p(\cdot)$ satisfying $1<p^-<p^+<+\infty$ and the so-called $\log$-H\"older regularity, that is
\begin{equation}\label{logholder_condition}
|p(x)-p(y)|\leq\frac{-C_0}{\log(|x-y|)}; \mbox{ for every }x,\;y\in\Omega \mbox{ with }|x-y|\leq\frac{1}{2},
\end{equation}
for some constant $C_0>0$.
\par In \cite[Theorem 8.2.4]{DHHR} the authors defined $W^{1,p(\cdot)}_0(\Omega)$ as the closure of Sobolev functions with compact support in $\Omega$ with respect to the norm in
$W^{1,p(\cdot)}(\Omega)$. They proved the Poincar\'e norm inequality for a regular bounded  domain and for exponent $p(\cdot)$ satisfying the following two conditions
\begin{equation}\label{logholder1}
\Big|\frac{1}{p(x)}-\frac{1}{p(y)}\Big|\leq\frac{C_1}{\log(e+\frac{1}{|x-y|})}
\end{equation}
and for some $p_\infty\in\mathbf{R}$
\begin{equation}\label{logholder2}
\Big|\frac{1}{p(x)}-p_\infty\Big|\leq\frac{C_2}{\log(e+|x|)}
\end{equation}
for every $x$, $y\in\Omega$, where $C_1>0$ and $C_2>0$ are constants.\\
Let us note in passing that the above two definitions of the Sobolev space $W^{1,p(\cdot)}_0(\Omega)$ coincide if $p(\cdot)$ is a measurable bounded exponent and if (\ref{logholder1}) and (\ref{logholder2}) are fulfilled (see \cite[Corollary 11.2.4]{DHHR}). Ciarlet and Dinca \cite{CD2011} proved the Poincar\'e norm inequality using an approach which does not rely on the density arguments.
\par In general in the variable exponent Sobolev spaces $W^{1,p(\cdot)}_0(\Omega)$ (defined as the norm closure of $C^\infty_0(\Omega)$ functions in $W^{1,p(\cdot)}(\Omega)$), the Poincar\'{e} integral inequality (\ref{poincare1}) with variable exponent $p(\cdot)$ instead of constant exponent $p$ fails to hold as it was shown in \cite[Example, pp.444-445]{FZ2001}. Indeed, if the variable exponent $p(\cdot)$ is a continuous function having a minimum or a maximum then an integral version of the Poincar\'{e} inequality can not be obtained (see \cite{FAN2005}). However, under a suitable monotony property on the variable exponent $p(\cdot)$ Maeda \cite{Maeda2008} proved the Poincar\'{e} integral inequality for $\mathcal{C}^1_0(\Omega)$-functions.
\par It is worth recalling that the Poincar\'{e} norm inequality in $W^{1,p(\cdot)}_0(\Omega)$ obtained in the aforementioned references requires the continuity of the variable exponent. Here we prove the Poincar\'{e} integral  inequality in Musielak spaces, and so the Poincar\'e norm inequality, by introducing some assumptions that don't require the continuity of the variable exponent when reducing to Sobolev spaces $W^{1,p(\cdot)}_0(\Omega)$.
\subsection{Structural assumptions}
In this subsection we give the definition of Musielak $\Phi$-functions and we introduce new systematic sufficient conditions which enable us to prove  Poincar\'e-type integral inequalities in Musielak spaces.
\begin{definition}\label{def:Phifn}($\phi$-function, $\Phi$-function). A real function $M$ : $\Omega\times\mathbf{R}^{+}\to\mathbf{R}^{+}$ is called a $\phi$-function, written $M\in\phi$, if $M(x,\cdot)$ is a nondecreasing and convex function for all $x\in\Omega$ with $M(x,0)=0$, $M(x,s)>0$ for $s>0$, $M(x,s)\rightarrow\infty \mbox{ as }s\rightarrow\infty$
	and $M(\cdot,s)$ is a measurable function for every $s\geq 0$.
	\par A $\phi$-function is called $\Phi$-function, written $M\in\Phi$, if furthermore it satisfies
	$$
	\lim_{s\to 0} \frac{M(x,s)}{s}=0\quad \mbox{ and }\quad \lim_{s\to \infty} \frac{M(x,s)}{s}=\infty.
	$$
\end{definition}
\par Throughout the paper, we consider $\Phi$-functions on which we assume at least one of the following fundamental regularity assumptions.
\begin{enumerate}[($\mathcal{M}$1)]
	\item \label{X1} There exists a function  $\varphi:\big[0, {1}/{2} ]\times\mathbf{R}^+\to\mathbf{R}^+$ such that $\varphi(\cdot,s)$ and $\varphi(x,\cdot)$ are nondecreasing functions and for all $x,y\in\overline{\Omega}$ with $|x-y|\leq\frac{1}{2}$ and for any constant $c>0$
	$$
	M(x,s)\leq\varphi(|x-y|,s)M(y,s),\quad\mbox{ with } \limsup_{\varepsilon\rightarrow0^+}\varphi(\varepsilon, c\varepsilon^{-N})<\infty.
	$$
	\item \label{monotonie}
	A $\Phi$-function $M$ is said to satisfy the $Y$-condition on a segment $[a,b]$ of the real line $\mathbf{R}$, if\\
	Either
	$$
	(Y_0) :
	\left\{
	\begin{array}{llr}
	\mbox{There exist } t_0\in\mathbf{R}^+\mbox{ and } 1\leq i\leq N \mbox{ such that the partial function }\\
	x_i\in [a,b]\mapsto M(x,t)\mbox{ changes constantly its monotony on both}\\
	\mbox{ sides of }t_0\; (\mbox{that is for }  t\geq t_0 \mbox{ and } t< t_0),
	\end{array}
	\right.
	$$
	Or
	$$
	(Y_\infty) :
	\left\{
	\begin{array}{lll}
	\mbox{There exists } 1\leq i\leq N
	\mbox{ such that for all }t\geq0, \mbox{ the partial function }\\
	x_i\in [a,b]\mapsto M(x,t) \mbox{ is monotone on }[a,b].
	\end{array}
	\right.
	$$
	Here, $x_i$ stands for the $i^{th}$ component of $x\in\Omega$.
\end{enumerate}
The highly challenging and important part of the analysis in Musielak spaces is giving a relevant structural condition yielding approximation properties of these nonstandard spaces. In general for a $\Phi$-function $M$, smooth functions are not dense in norm in the Musielak space $W^mL_M(\Omega)$. The authors \cite{AYGS2017} introduced the condition \ref{X1} to study the problem of density of smooth functions in Musielak spaces and they showed that this condition unify and improve the known results in Orlicz-Sobolev spaces as well as the variable exponent Sobolev spaces. In fact, the condition \ref{X1} holds trivially in the case of Orlicz spaces while in the case of variable exponent Sobolev spaces \ref{X1} holds if we choose
$$
\varphi(\tau,s)= \max\left\{s^{\sigma(\tau)},s^{-\sigma(\tau)}\right\}.
$$
When $\sigma(\tau)=-{c}/{\log\tau}$, with $0<\tau\leq {1}/{2}$, we obtain the $\log$-H\"{o}lder continuity condition (\ref{logholder_condition}). Nonetheless, we can choose various $\varphi$s. For more examples of $\Phi$-functions satisfying \ref{X1} we refer to \cite{AYGS2017}.
\begin{remark}\label{remark1}
	\begin{enumerate}
		\item In the case where $M(x,t)=t^{p(x)}$, the assumption $(Y_0)$ prevents the variable exponent $p(\cdot)$ to get a local extremum while $(Y_\infty)$ is not satisfied unless $p(\cdot)$ is a constant function.
		\item Let us consider the double phase function $M(x,t)=t^{p}+a(x)t^q$. If there is $1\leq i\leq N$ such that the function $x_i \mapsto a(x)$ is monotone then $M$ satisfies obviously $(Y_\infty)$ and so \ref{monotonie}. If $x\mapsto a(x)$ is not a constant function then the double phase function $M$ can not satisfy $(Y_0)$.
		\item If $1<p(\cdot)<+\infty$ and there exists $1\leq i\leq N$ such that the function $x_i\mapsto p(x)$ is monotone on a compact subset of the real line $\mathbf{R}$, then the following $\Phi$-functions
		$$
		M_1(x,t)=t^{p(x)},\qquad M_2(x,t)=t^{p(x)}\log (e+t), \qquad M_3(x,t)=e^{t^{p(x)}}-1,
		$$
		satisfy \ref{monotonie}.
	\end{enumerate}
\end{remark}
We note in passing here, that the assumption \ref{monotonie} covers the one given in \cite{Maeda2008}.\\
In what follows, we will use the fact that any $\Phi$ function $M$ is locally integrable in $\Omega$, that is to say for any constant number $c>0$ and for every compact set $K\subset\Omega$
\begin{equation}\label{incBM}
\int_K M(x,c)dx<\infty.
\end{equation}
We note that (\ref{incBM}) is obviously satisfied. Indeed, defining the increasing sequence $\{\Omega_j\}_{j\geq1}$
$$
\Omega_j=\bigg\{x\in\Omega: |x|< j, dist(x,\Omega^c)>\frac{1}{j}\bigg\},
$$
one has $\Omega=\cup_{j=1}^\infty \Omega_j$. So that for any compact subset $K\subset\Omega$ there is a finite recovering such that $K\subset\cup_{j=1}^p \Omega_j$. Then
$$
\int_{K}M(x,c)dx\leq\sum_{j=1}^p\int_{\Omega_j}M(x,c)dx.
$$
Therefore, (\ref{incBM}) follows from \cite[p. 64]{Kaminska1985}.
\subsection{Main results}
In this subsection we give our main results. Let $M$ be a $\Phi$-function and let $M^\ast$ be its complementary $\Phi$-function (see (\ref{complementaryfunction}) hereafter). Define the space $W^{m}_0L_M(\Omega)$ to be the $\sigma(\Pi L_M, \Pi E_{M^{\ast}})$ closure of $\mathcal{C}^{\infty}_{0}(\Omega)$ in $W^{m}L_M(\Omega)$.
The first result we obtain concerns Poincar\'e-type  inequalities in the Musielak spaces $W^m_0L_M(\Omega)$.
\begin{theorem}\label{thpoincare}
	Let $\Omega$ be a bounded open subset in $\mathbf{R}^N$ having the segment property and let $M\in\Phi$ satisfies \ref{X1} and  \ref{monotonie}. Then there exists a constant $c_{m,\Omega}$ depending only on $m$ and $\Omega$ such that for every $u\in W^m_0L_M(\Omega)$
	\begin{equation}\label{poincaremod}
	\int_\Omega \sum_{|\alpha|< m} M(x,|D^\alpha u|)dx\leq \int_\Omega \sum_{|\alpha|= m} M(x,c_{m,\Omega}|D^\alpha u|)dx.
	\end{equation}
	Moreover, for every $u\in W^m_0L_M(\Omega)$
	\begin{equation}\label{poincarenorm}
	\sum_{|\alpha|< m} \|D^\alpha u\|_{M,\Omega}\leq C(m,\Omega)\sum_{|\alpha|= m} \|D^\alpha u\|_{M,\Omega},
	\end{equation}
	where $C(m,\Omega)$ is a constant depending only on $m$ and $\Omega$.
\end{theorem}
\par In the framework of Orlicz spaces, Theorem \ref{thpoincare} was proved by Gossez \cite[Lemma 5.7]{GJP1} where only the definition of the space $W^{m}_0L_{\varphi}$, defined as the closure of $\mathcal{C}^\infty_0(\Omega)$-functions with respect to the weak-$\ast$ topology $\sigma(\Pi L_\varphi,\Pi E_{\varphi^\ast})$, is used to get the Poincar\'e integral inequality without assuming the segment property on the bounded open $\Omega$. As in the classical way, the Poincar\'e integral inequality was first proved for smooth functions and then (\ref{poincaremod}) follows by a density argument based on mollifications.\\
In general, the shift operator is not acting between Musielak spaces (see \cite[example 2.9 and Theorem 2.10]{KR}). So we face a major difficulty in using mollification and then we can not use the same approach as in \cite{GJP1}.
\par Our contribution to overcame this problem consists in using the regularity condition \ref{X1} on the $\Phi$-function $M$ and the segment property on the domain $\Omega$ (see Definition \ref{segmpropdef}). Those conditions enable us to get the modular density of $C^\infty_0$-functions in $W^m_0L_M(\Omega)$
(see \cite{AYGS2017})  which we use in the proof of Theorem \ref{thpoincare}
instead the weak-$\ast$ density as it was done in \cite{GJP1}.
\begin{remark}
	A direct consequence of the inequality (\ref{poincarenorm}), is that the following two norms $\|D^\alpha u\|_{m,M,\Omega}$ and $\|D^m u\|_{M,\Omega}$ are equivalent on $W^{m}_{0}L_M(\Omega)$.
\end{remark}
The following theorem concerns the Poincar\'e integral inequality in the Musielak-Sobolev space $W^m_0E_M(\Omega)$ defined as the norm closure of $\mathcal{C}^{\infty}_{0}$-functions in $W^mE_M(\Omega)$.
\begin{theorem}\label{cor.poincare}
	Let $\Omega$ be a bounded open subset in $\mathbf{R}^N$ and let $M\in\Phi$ satisfies  \ref{monotonie}. The inequality (\ref{poincaremod}) holds true for every $ u\in W^m_0E_M(\Omega)$ and then so is (\ref{poincarenorm}).
\end{theorem}
We note here that by the definition of $W^m_0E_M(\Omega)$, we do not need to assume in the above Theorem \ref{cor.poincare} the segment property on $\Omega$ and the condition \ref{X1} on the $\Phi$-function $M$.
Therefore, in view of Remark \ref{remark1} the result we obtain covers the Poincar\'e integral inequality obtained by Maeda \cite{Maeda2008} for $\mathcal{C}^{1}_{0}$-functions in the case where the variable exponent $p(\cdot)$ is assumed to satisfy a monotony condition.
\par Let $K^{m}_0L_M(\Omega)$ be the norm closure of the set of $W^{m}L_M(\Omega)$ functions with compact support in $\Omega$. In the particular case  $M(x,t)=t^{p(x)}$, $K^{m}_0L_M(\Omega)$ is nothing but the space $W^{m,p(x)}_0(\Omega)$ defined in \cite{DHHR}.
\begin{theorem}\label{thzeroboundary}
	Let $\Omega$ be an open subset in $\mathbf{R}^N$ and let $M\in\Phi$ satisfies  \ref{X1}. Then $K^{m}_0L_M(\Omega)$ coincides with $W^{m}_0E_M(\Omega)$. Furthermore, if $\Omega$ is bounded and $M$ satisfies \ref{monotonie} then (\ref{poincaremod}) and (\ref{poincarenorm}) are fulfilled.
\end{theorem}
Now, the remaining question is how to provide a satisfactory generalization of the Poincar\'e inequality for constant exponent, because the equality (\ref{tr}) is substituted by the inclusion
$$
W^{1}_{0}L_M(\Omega)=\overline{\mathcal{C}^\infty_0(\Omega)}^{\sigma(\Pi L_M,\Pi E_{M^\ast})}\subset\big\{u\in W^{1}L_M(\Omega) : tr(u)=0 \mbox{ on } \partial\Omega\big\}
$$
which may be strict in general unless additional conditions are imposed on the $\Phi$-function $M$. Note here that for every $u\in W^{1}L_{M}(\Omega)$ the trace $tr (u)=u|_{\partial\Omega}$ is well defined. Indeed, if $\Omega$ is of finite Lebesgue measure one has
$W^{1}L_{M}(\Omega)\hookrightarrow W^{1,1}(\Omega)$ and by the Gagliardo trace theorem (see \cite{gag}) we have the embedding
$W^{1,1}(\Omega)\hookrightarrow L^{1}(\partial\Omega)$. Hence, we conclude that for all $u\in W^{1}L_{M}(\Omega)$ there holds $u|_{\partial\Omega}\in L^{1}(\partial\Omega)$. We give the answer in the following theorem.
\begin{theorem}\label{th5.8}
	Let $\Omega$ be a bounded open subset in $\mathbf{R}^N$ having the segment property. Assume that $M\in\Phi$ satisfies \ref{X1}. Then, we get
	$$
	W^{m}_{0}L_M(\Omega)=W^{m,1}_{0}(\Omega)\cap W^{m}L_M(\Omega).
	$$
	If furthermore $\Omega$ has a Lipschitz boundary $\partial\Omega$, then we obtain
	\begin{equation}\label{th5.8_eq2}
	W^{1}_{0}L_M(\Omega)=\{u\in W^{1}L_M(\Omega) : \mbox{ tr}(u)=0 \mbox{ on } \partial\Omega\}.
	\end{equation}
\end{theorem}
\subsection{Organization of the paper }
In section \ref{section2} we review some basic facts we use about Musielak spaces. Further details can be found in the standard monograph by J. Musielak 
\cite{MJ} and the papers by Kami\'nska \cite{Kaminska1981,Kaminska1984,Kaminska1985}. Section~\ref{section4} is devoted to the proof of the main results.
\section{Musielak Structure}\label{section2}
In the following section we give a brief basic review on Musielak-Orlicz spaces.
For $M\in\phi$, the Musielak-Orlicz space $L_M(\Omega)$ (resp. $E_M(\Omega)$) is defined as the set of all measurable functions $u:\Omega\rightarrow \mathbf{R}$ such that $\int_\Omega M(x,|u(x)|/\lambda)dx<+\infty$ for some $\lambda>0$ ($\mbox{resp. for all } \lambda>0$). Equipped with the Luxemburg norm
$$
\|u\|_{M,\Omega}=\inf\bigg\{\lambda>0: \int_\Omega M\bigg(x,\frac{|u(x)|}{\lambda}\bigg)dx\leq 1\bigg\}.
$$
$L_M(\Omega)$ is a Banach space \cite[Theorem 7.7]{MJ} and $E_M(\Omega)$ is its closed subset. Define $M^{\ast}: \Omega\times\mathbf{R}^{+}\to\mathbf{R}^{+}$ by
\begin{equation}\label{complementaryfunction}
M^{\ast}(x,s)=\sup_{t\geq 0}\{st-M(x,t)\} \mbox{ for all } s\geq 0\mbox{ and all } x\in \Omega.
\end{equation}
$M^{\ast}$ is also a $\Phi$-function and is called the complementary function to $M$ in the sense of Young. Moreover, we have the following Young inequality
$$
uv\leq M(x,u)+ M^{\ast}(x,v),\quad \forall u,v\geq 0, \forall x\in\Omega,
$$
from which we easily get the H\"{o}lder inequality
$$
\int_\Omega |uv| dx\leq2\|u\|_{M,\Omega}\|v\|_{M^{\ast},\Omega}
$$
for all $u\in L_M(\Omega)$ and $v\in L_{M^{\ast}}(\Omega)$.
We say that $\{u_k\}_k$ converges to $u$ in norm in $L_M(\Omega)$, if
$\|u_k-u\|_{M,\Omega}\rightarrow 0$ as $k\rightarrow \infty$. The notion of the modular convergence is given in the following definition.
\begin{definition}[Modular convergence]\label{def:mod-conv}
	A sequence $\{u_k\}_k$ is said to converge modularly to $u$ in $L_M(\Omega)$
	if there exists $\lambda>0$ such that
	$$
	\rho_{M}((u_k-u)/\lambda):=\int_{\Omega}M\left(x,{|u_k-u|}/{\lambda}\right)\, dx\rightarrow 0 \mbox{ as } k\rightarrow \infty.
	$$
\end{definition}
For a positive integer $m$, we define the Musielak-Orlicz-Sobolev spaces  $W^{m}L_M(\Omega)$ and $W^{m}E_M(\Omega)$ as follows
$$
W^{m}L_M(\Omega)=\Big\{u\in L_M(\Omega): D^\alpha u\in L_M(\Omega), |\alpha|\leq m\Big\},
$$
$$
W^{m}E_M(\Omega)=\Big\{u\in E_M(\Omega): D^\alpha u \in E_M(\Omega), |\alpha|\leq m\Big\},
$$
where $\alpha=(\alpha_1,\alpha_2,\cdots,\alpha_N)$, $|\alpha|=|\alpha_1|+|\alpha_2|+\cdots+|\alpha_N|$ and $D^\alpha =\frac{\partial^{|\alpha|}}{\partial^{\alpha_1}_{x_1}\cdots\partial^{\alpha_N}_{x_N}}$ stands for the distributional derivatives. Observe that by (\ref{incBM}) the function $x\to M^\ast(x,c)$ always belongs to $L_{loc}^1(\Omega)$ for every constant number $c\geq0$. Hence, for an arbitrary compact subset $K$ of $\Omega$ we can write by H\"older's inequality in Musielak spaces
$$
\int_{K}|u(x)|dx\leq 2\|u\|_{0,M,\Omega}\|\chi_{K}\|_{0,M^\ast,\Omega}
\leq 2\Big(\int_{K}M^\ast(x,1)dx+1\Big)\|u\|_{0,M,\Omega}
$$
which yields $L_M(\Omega)\subset L_{loc}^1(\Omega)$. Therefore, endowed with the Luxemburg norm
$$
\|u\|_{m,M,\Omega}=\inf\bigg\{\lambda>0 : \sum_{ |\alpha|\leq m} \rho_{M}(D^\alpha u/\lambda) \leq 1\bigg\}.
$$
$\big(W^{m}L_M(\Omega), \|u\|_{m,M,\Omega}\big)$ is a Banach space.
We will always identify the space $W^mL_M(\Omega)$ to a subspace of the product $\Pi_{|\alpha|\leq m} L_M=\Pi L_M$.
\begin{definition}[Segment property]\label{segmpropdef}
	A domain $\Omega$ is said to satisfy the  {segment property}, if there exist a finite open covering $\{\theta\}_{i=1}^k$ of $\overline{\Omega}$ and a corresponding non-zero vectors $y_i\in \mathbf{R}^N$ such that  $(\overline{\Omega}\cap\theta_{i})+ty_{i}\subset\Omega$ for all $t\in (0,1)$ and $i=1,\dots,k$.
\end{definition}
This condition holds, for example, if $\Omega$ is a bounded Lipschitz domain (cf. \cite{AF}). By convention, the empty set satisfies the segment property.
\section{Proof of main results}\label{section4}
\subsubsection*{Proof of Theorem \ref{thpoincare}}
Let $d$ be the diameter of $\Omega$. As \ref{monotonie} is concerned and without loss of generality, we can assume that $i=1$. Being $\Omega$ bounded, using a translation if necessary, we may assume that it is contained in the strip
$\Omega\subset\{(x_1,x^{\prime})\in[0,d]\times\mathbf{R}^{N-1}\}$. Let $\partial_{1}:=\frac{\partial }{\partial x_1}$ stands for the partial derivative operator with respect to $x_1$ and let us first assume that $u\in \mathcal{C}^\infty_0(\Omega)$.
\par\noindent\textbf{Part 1 : } We assume that there exists $t_0\in\mathbf{R}^{+}$ such that the function $x_1\in [0,d]\mapsto M\Big((x_1,x^{\prime}),t\Big)$ changes the variation on both sides of $t_0$.\\
\par\noindent\textbf{Case 1. } Assume that $x_1\in [0,d]\mapsto M\Big((x_1,x^{\prime}),t\Big)$
is non-decreasing for $t\leq t_0$ and non-increasing for $t_0<t$.\\
Defining the two sets
$$
E_1=\Big\{\xi\in [0,d]:|\partial_{1}u(\xi,x^\prime)|\leq \frac{1}{d}t_0\Big\}
\mbox{ and }
E_2=\Big\{\xi\in [0,d]:|\partial_{1}u(\xi,x^\prime)|> \frac{1}{d}t_0\Big\},
$$
we can write
$$
\begin{array}{lll}
u(x_1,x^\prime)&=&u(x_1,x^\prime)\chi_{E_1}(x_1)+u(x_1,x^\prime)\chi_{E_2}(x_1)\\
&=&-\int_{x_1}^d \partial_{1}\Big(u(\xi,x^\prime)\chi_{E_1}(\xi)\Big)d\xi+\int_{0}^{x_1} \partial_{1}\Big(u(\xi,x^\prime)\chi_{E_2}(\xi)\Big)d\xi\\
&=&-\int_{[x_1,d]\cap E_1} \partial_{1}u(\xi,x^\prime)d\xi+\int_{[0,x_1]\cap E_2}\partial_{1}u(\xi,x^\prime)d\xi.
\end{array}
$$
Thus,
$$
|u(x_1,x^\prime)|\leq\int_{[x_1,d]\cap E_1} \Big|\partial_{1}u(\xi,x^\prime)\Big|d\xi+\int_{[0,x_1]\cap E_2} \Big|\partial_{1}u(\xi,x^\prime)\Big|d\xi.
$$
Then, the convexity of the $\Phi$-function $M$ and Jensen's inequality enable us to write
$$
\begin{array}{lll}
M\Big(x,|u(x_1,x^\prime)|\Big)&\leq& \frac{1}{2d}\int_0^d M\Big((x_1,x^\prime),2d\Big|\partial_{1}u(\xi,x^\prime)\Big|\chi_{[x_1,d]\cap E_1}(\xi)\Big)d\xi\\
&+&\frac{1}{2d}\int_0^d M\Big((x_1,x^\prime),2d\Big|\partial_{1}u(\xi,x^\prime)\Big|\chi_{[0,x_1]\cap E_2}(\xi)\Big)d\xi\\
&\leq& \frac{1}{2d}\int_{[x_1,d]\cap E_1} M\Big((\xi,x^\prime),2d\Big|\partial_{1}u(\xi,x^\prime)\Big|\Big)d\xi\\
&+&\frac{1}{2d}\int_{[0,x_1]\cap E_2} M\Big((\xi,x^\prime),2d\Big|\partial_{1}u(\xi,x^\prime)\Big|\Big)d\xi\\
&\leq& \frac{1}{2d}\int_{0}^d M\Big((\xi,x^\prime),2d\Big|\partial_{1}u(\xi,x^\prime)\Big|\Big)d\xi.
\end{array}
$$
Integrating successively with respect to $x^\prime$ and $x_1$, we obtain
\begin{equation}\label{poinc_part1case1}
\int_\Omega M\Big(x,|u(x)|\Big)dx\leq \frac{1}{2}\int_\Omega M\Big(x,2d\Big|\partial_{1}u(x)\Big|\Big)dx.
\end{equation}
\par\noindent\textbf{Case 2. } Assume that $x_1\in [0,d]\mapsto M\Big((x_1,x^{\prime}),t\Big)$ is non-increasing on $t\leq t_0$ and non-decreasing on $t_0<t$.\\
We can write
$$
\begin{array}{lll}
u(x_1,x^\prime)&=&u(x_1,x^\prime)\chi_{E_1}(x_1)+u(x_1,x^\prime)\chi_{E_2}(x_1)\\
&=&\int_{0}^{x_1}  \partial_{1}\Big(u(\xi,x^\prime)\chi_{E_1}(\xi)\Big)d\xi-\int_{x_1}^d \partial_{1}\Big(u(\xi,x^\prime)\chi_{E_2}(\xi)\Big)d\xi\\
&=&\int_{[0,x_1]\cap E_1} \partial_{1}u(\xi,x^\prime)d\xi-\int_{[x_1,d]\cap E_2} \partial_{1}u(\xi,x^\prime)d\xi,
\end{array}
$$
which implies
$$
|u(x_1,x^\prime)|\leq\int_{[0,x_1]\cap E_1} \Big|\partial_{1}u(\xi,x^\prime)\Big|d\xi+\int_{[x_1,d]\cap E_2} \Big|\partial_{1}u(\xi,x^\prime)\Big|d\xi.
$$
Once again the convexity of the $\Phi$-function $M$ and Jensen's inequality enable us to write
$$
\begin{array}{lll}
M\Big(x,|u(x_1,x^\prime)|\Big)&\leq& \frac{1}{2d}\int_0^d M\Big((x_1,x^\prime),2d\Big|\partial_{1}u (\xi,x^\prime)\Big|\chi_{[0,x_1]\cap E_1}(\xi)\Big)d\xi\\
&+&\frac{1}{2d}\int_0^d M\Big((x_1,x^\prime),2d\Big|\partial_{1}u(\xi,x^\prime)\Big|\chi_{[x_1,d]\cap E_2}(\xi)\Big)d\xi\\
&\leq& \frac{1}{2d}\int_{[0,x_1]\cap E_1} M\Big((\xi,x^\prime),2d\Big|\partial_{1}u(\xi,x^\prime)\Big|\Big)d\xi\\
&+&\frac{1}{2d}\int_{[x_1,d]\cap E_2} M\Big((\xi,x^\prime),2d\Big|\partial_{1}u(\xi,x^\prime)\Big|\Big)d\xi\\
&\leq& \frac{1}{2d}\int_{0}^d M\Big((\xi,x^\prime),2d\Big|\partial_{1}u(\xi,x^\prime)\Big|\Big)d\xi.
\end{array}
$$
Integrating successively with respect to $x^\prime$ and $x_1$, we obtain
$$
\int_\Omega M\Big(x,|u(x)|\Big)dx\leq \frac{1}{2}\int_\Omega M\Big(x,2d\Big|\partial_{1}u(x)\Big|\Big)dx.
$$
\par\noindent\textbf{Part 2 :} Assume now that for all $t\geq0$, the function $x_1\in [0,d]\mapsto M\Big((x_1,x^{\prime}),t\Big)$ is monotone.
\par\noindent\textbf{Case 1. } Assume first that $x_1\in [0,d]\mapsto M\Big((x_1,x^{\prime}),t\Big)$ is non-increasing.\\
By using Jensen's inequality we get
$$
\begin{array}{lcl}
M\Big(x,|u(x_1,x^\prime)|\Big) &\leq& M\Big((x_1,x^{\prime}),\int_{0}^{x_1}\Big|\partial_{1}u(\xi,x^{\prime})\Big|d\xi\Big)\\
&\leq& M\Big((x_1,x^{\prime}),\int_{0}^{d}\Big|\partial_{1}u(\xi,x^{\prime})\Big|\chi_{[0,x_1]}(\xi)d\xi\Big)\\
&\leq&\frac{1}{d}\int_{0}^{d}M\Big((x_1,x^{\prime}), d\Big|\partial_{1}u(\xi,x^{\prime})\Big|\chi_{[0,x_1]}(\xi)\Big)d\xi\\
&\leq&\frac{1}{d} \int_{0}^{x_1}M\Big((\xi,x^{\prime}),d\Big|\partial_{1}u(\xi,x^{\prime})\Big|\Big)d\xi\\
&\leq&\frac{1}{d} \int_{0}^{d}M\Big((\xi,x^{\prime}),d\Big|\partial_{1}u(\xi,x^{\prime})\Big|\Big)d\xi.
\end{array}
$$
Integrating successively with respect to $x^\prime$ and $x_1$, we obtain
$$
\int_\Omega M\Big(x,|u(x)|\Big)dx\leq \int_\Omega M\Big(x,d\Big|\partial_{1}u(x)\Big|\Big)dx.
$$
\par\textbf{Case 2. } Assume that $x_1\in [0,d]\mapsto M\Big((x_1,x^{\prime}),t\Big)$ is non-decreasing.\\
By virtue of Jensen's inequality we can write
$$
\begin{array}{lcl}
M\Big(x,|u(x_1,x^\prime)|\Big) &\leq& M\Big((x_1,x^{\prime}),\int_{x_1}^{d}\Big|\partial_{1}u(\xi,x^{\prime})\Big|d\xi\Big)\\
&\leq& M\Big((x_1,x^{\prime}),\int_{0}^{d}\Big|\partial_{1}u(\xi,x^{\prime})\Big|\chi_{[x_1,d]}(\xi)d\xi\Big)\\
&\leq&\frac{1}{d} \int_{0}^{d}M\Big((x_1,x^{\prime}),d\Big|\partial_{1}u(\xi,x^{\prime})\Big|\chi_{[x_1,d]}(\xi)\Big)d\xi\\
&\leq&\frac{1}{d} \int_{x_1}^{d}M\Big((\xi,x^{\prime}),d\Big|\partial_{1}u(\xi,x^{\prime})\Big|\Big)d\xi\\
&\leq&\frac{1}{d} \int_{0}^{d}M\Big((\xi,x^{\prime}),d|\partial_{1}u(\xi,x^{\prime})|\Big)d\xi.
\end{array}
$$
Integrating successively with respect to $x^\prime$ and $x_1$, we obtain
\begin{equation}\label{poinc_part2case2}
\int_\Omega M\Big(x,|u(x)|\Big)dx\leq \int_\Omega M\Big(x,d\Big|\partial_{1}u(x)\Big|\Big)dx.
\end{equation}
\par To sum up, from (\ref{poinc_part1case1})-(\ref{poinc_part2case2}), we obtain
\begin{equation}\label{eq2.1.1}
\int_{\Omega}M\Big(x,|u(x)|\Big)dx\leq \int_{\Omega}M\Big(x,2d\Big|\partial_{1}u(x)\Big|\Big)dx,
\end{equation}
for all $u\in \mathcal{C}^\infty_0(\Omega)$.
\par\noindent Let now $u\in W^{1}_0L_M(\Omega)$ be arbitrary. By \cite[Theorem 3]{AYGS2017} there exist $\lambda>0$ and a sequence of functions $u_k\in\mathcal{C}^\infty_0(\Omega)$ such that
$$
\int_\Omega M\Big(x,\frac{|u_k(x)-u(x)|}{\lambda}\Big)dx+\int_\Omega M\Big(x,\frac{|\nabla u_k(x)-\nabla u(x)|}{\lambda}\Big)dx\rightarrow 0
$$
as $k\to+\infty$. Hence, up to a subsequence still again indexed by $k$, we can assume that $u_k\rightarrow u$ a.e. in $\Omega$. Then, using (\ref{eq2.1.1}) we can write
$$
\begin{array}{lcl}
\int_\Omega M\Big(x, \frac{|u(x)|}{4\lambda d}\Big)dx &\leq&\lim_{
	k\rightarrow+\infty} \inf\int_\Omega M\Big(x, \frac{|u_{k}(x)|}{4\lambda d}\Big)dx\\
&\leq&\lim_{k\rightarrow+\infty} \inf\int_\Omega M\Big(x, \frac{1}{2\lambda}\Big|\partial_1 u_k(x)\Big|\Big)dx\\
&\leq& \frac{1}{2}\lim_{k\rightarrow+\infty} \inf\int_\Omega M\Big(x, \frac{1}{\lambda}\Big|\partial_1 u_k(x)-\partial_1u(x)\Big|\Big)dx\\
&&+\frac{1}{2}\int_\Omega M\Big(x,\frac{1}{\lambda}\Big|\partial_1 u(x)\Big|\Big)dx\\
&\leq& \frac{1}{2}\int_\Omega M\Big(x,\frac{1}{\lambda}\Big|\partial_1 u(x)\Big|\Big)dx.
\end{array}
$$
Thus, (\ref{poincaremod}) is proved.
Let us now prove the inequality (\ref{poincarenorm}). For $u\in W^m_0L_M(\Omega)$, it can be checked easily from (\ref{poincaremod}) that
$$
\sum_{|\alpha|<m}\int_\Omega M\Big(x,\frac{|D^{\alpha} u(x)|}{C(m,\Omega)
	{\sum_{|\beta|=m}}\|D^{\beta} u\|_{M,\Omega}}\Big)dx
\leq 1,
$$
where $C(m,\Omega)=c_{m,\Omega}\Big(1+\sum_{|\beta|=m}1\Big)$ depending only on $m$ and $\Omega$. The proof of Theorem \ref{thpoincare} is then achieved.\qed
\subsection*{Proof of Theorem \ref{cor.poincare}}
Let $u\in W^m_0E_M(\Omega)$. By the definition of $W^m_0E_M(\Omega)$, there exists a sequence $\{u_k\}_k$ of $\mathcal{C}^{\infty}_0(\Omega)$ functions such that $D^\alpha u_k\to D^\alpha u$ for all $|\alpha|\leq m$ with respect to the norm topology in $L_M(\Omega)$. As the norm convergence implies the modular one, one has $D^\alpha u_k\to D^\alpha u$ for all $|\alpha|\leq m$ with respect to the modular topology. Therefore, we get the result by following exactly the same lines of the proof of Theorem \ref{thpoincare}.\qed
\subsection*{Proof of Theorem \ref{thzeroboundary}}
The embedding $W^m_0E_M(\Omega)\subset K^m_0L_M(\Omega)$ is obviously satisfied. It only remains to show  that $K^m_0L_M(\Omega)\subset  W^m_0E_M(\Omega)$ holds true. Let $u\in K^m_0L_M(\Omega)$ and let $\eta>0$ be arbitrary.  We will show that there is a sequence $v\in \mathcal{C}^\infty_0(\Omega)$ such
\begin{equation}\label{cv}
\sum_{|\alpha|\leq m}\|D^\alpha u-D^\alpha v\|_{M,\Omega}\leq\eta.
\end{equation}
By the definition of $K^m_0L_M(\Omega)$ there exist a sequence $\{u_k\}$ in $W^{m}L_M(\Omega)$ of compactly supported functions in $\Omega$ and $k_{\alpha}>0$ such that for all $k>k_{\alpha}$ and $|\alpha|\leq m$ we have
$$
\|D^\alpha u-D^\alpha u_k\|_{M,\Omega}\leq\frac{\eta}{2K}.
$$
where $K$ is  the total number of multi-indices with $|\alpha|\leq m$. Now by using \cite[Lemma 12]{AYGS2017} there exist a sequence $\{u^n_k\}$ in $\mathcal{C}^\infty_0(\Omega)$ and $n^{k,\alpha}_{\eta}>0$ such that for all $n>n^{k,\alpha}_{\eta}$ and $|\alpha|\leq m$ we have
$$
\|D^\alpha u_k-D^\alpha u^n_k\|_{M,\Omega}\leq\frac{\eta}{2K}.
$$
By the triangle inequality we get for all $k>k^{\alpha}_0$ and $n>n^{k,\alpha}_{\eta}$
$$
\begin{array}{lll}
\|D^\alpha u-D^\alpha u^n_k\|_{M,\Omega}&\leq&
\|D^\alpha u-D^\alpha u_k)\|_{M,\Omega}+ \|D^\alpha u_k-D^\alpha u^n_k\|_{M,\Omega}\\
&\leq&\frac{\eta}{K},
\end{array}
$$
Hence follows (\ref{cv}) and then the following identification
$$
K^{m}_0L_M(\Omega)=W^{m}_0E_M(\Omega).
$$
Therefore, Theorem \ref{thzeroboundary} follows immediately from Theorem \ref{cor.poincare}.\qed
\subsubsection*{Proof of Theorem \ref{th5.8}}
We begin first by showing that $W^{m}_{0}L_M(\Omega)\subset W^{m,1}_{0}(\Omega)\cap W^{m}L_M(\Omega)$. Since $W^{m}_{0}L_M(\Omega)$ is a subset of $W^{m}L_M(\Omega)$, it's sufficient to check that for all function $u$ belonging to $W^{m}_{0}L_M(\Omega)$ we have $u\in W^{m,1}_{0}(\Omega)$.
Let $u\in W^{m}_{0}L_M(\Omega)$. By \cite[Theorem 3]{AYGS2017} there exist $\lambda>0$ and a sequence $u_k\in \mathcal{C}^{\infty}_{0}(\Omega)$ such that
$$
\sum_{|\alpha|\leq m}\int_{\Omega}M (x,|D^\alpha u_k(x)-D^\alpha u(x)|/\lambda )dx\rightarrow 0 \mbox{ as } k\rightarrow \infty.
$$
Therefore, for all $|\alpha|\leq m$
$$
\int_{\Omega}M (x,|D^\alpha u_k(x)-D^\alpha u(x)|/\lambda )dx\rightarrow 0 \mbox{ as } k\rightarrow \infty.
$$
Thus, for a subsequence still denoted by $u_k$, we can assume
$$
D^\alpha u_k\to D^\alpha u \mbox{ a.e. in }\Omega.
$$
Applying Vitali's theorem we obtain
$$
\int_{\Omega}|D^\alpha u_k(x)-D^\alpha u(x)|dx \rightarrow 0 \mbox{ as } k\rightarrow \infty,
$$
which implies that $u\in W^{m,1}_0(\Omega)$.
\par Conversely, we should prove that $W^{m,1}_{0}(\Omega)\cap W^{m}L_M(\Omega)  \subset W^{m}_{0}L_M(\Omega)$.
We will show that for $u\in W^{m,1}_{0}(\Omega)\cap W^{m}L_M(\Omega)$ there exist a sequence
$v\in \mathcal{C}^\infty_0(\Omega)$ such that $v$ converges in the modular sense to $u$ in $W^mL_M(\Omega)$ and then we conclude by using \cite[Theorem 3]{AYGS2017}.
Let us denote by $\widetilde{u}$ the extension of $u$ by zero outside $\Omega$. Since $u$ belongs to  $W^{m,1}_{0}(\Omega)$ it yields,
$\widetilde{u}\in W^{m,1}(\mathbf{R}^N)$ and $\widetilde{D^\alpha u}=D^\alpha\widetilde{u}$ in the distributional sense and a.e. in
$\mathbf{R}^N$ (see \cite[Lemma 3.27]{AF}) and so $\widetilde{u}\in W^{m}L_M(\mathbf{R}^N)$. Then by \cite[Lemma 3]{AYGS2017} we can assume
that $u$ has compact support $K\subset\overline{\Omega}$. We will distinguish  the two cases: either $K\subset\Omega$ or $K\cap\partial\Omega\neq\emptyset$.
If $K\subset\Omega$ then we get the desired inclusion by \cite[Lemma 12]{AYGS2017}. If $K\cap\partial\Omega\neq\emptyset$, then, as in the proof
of \cite[Theorem 2]{AYGS2017}, there exist a finite collection $\{\widehat{\theta}_i\}_{i=1}^{k}$ covering the compact set $K\cap\partial\Omega$
and an open covering  $\{\theta^{\prime}_{i}\}^{k}_{i=0}$ of $K$ with $\theta^{\prime}_{i}$ has a compact closure in $\widehat{\theta}_i$
for $i=0,1,\cdots,k$. Then $u$ can be splitted into finitely-many pieces $u_i$,  such that $u=\sum_{i=1}^k u_i$ with
$\supp\, u_{i}\subset \theta_{i}^{\prime}$, $i=0,1,\cdots,k$.\\
For $i=0$, we consider $\supp\, u_{0}\subset\theta^{\prime}_0\subset \Omega$, then as for the first case by \cite[Lemma 12]{AYGS2017} there exist
$\varepsilon_0>0$ small enough ($\varepsilon_0<dist(\theta^{\prime}_0,\partial\Omega)$), such the regularized function $v_0=J_{\varepsilon_0}*u_{0}$
belongs to $\mathcal{C}^{\infty}_{0}(\Omega)$ and converge in modular since to $u$ in $W^mL_M(\Omega)$.\\
For $1\leq i\leq k$ fix. Let $z_i$ be a non-zero vector associated to $\widehat{\theta}_i$ by the segment property and let $r_i\in(0,1)$ be such that
$$
0<r_i<\min\{1/(|z_i|+1),dist(\theta_{i}^{\prime},\partial\widehat{\theta}_i)|z_i|^{-1}\}.
$$
Define \[(u_{i})_{-r_i}(x)=u_{i}(x-r_iz_i).\] and choose
$$
\varepsilon_i<dist\bigg((\theta^\prime_i\cap \overline{\Omega})+r_iz_i, \mathbf{R}^N\setminus\Omega\bigg).
$$
We define then the sequences
$$
v_i^{\varepsilon_i,r_i}(x)=J_{\varepsilon_i}*(u_{i})_{-r_i}=\int_{B(0,1)} J(y)u_i(x-r_iz_i-\varepsilon_i y)dy
$$
and
$$
v(x)=\sum_{i=1}^{k}v_i^{\varepsilon_i,r_i}(x)+J_{\varepsilon_0}\ast u_0(x),
$$
therefore,  $v\in \mathcal{C}^\infty_0(\Omega)$. Then, arguing similarly as in the proof of \cite[Theorem 2]{AYGS2017}, we prove that $v$ converges to
$u$ in $W^mL_M(\Omega)$ with respect to the modular convergence. This implies that $u$ belongs to  $W^m_0L_M(\Omega)$.\\
To check (\ref{th5.8_eq2}) observe that $\{u\in W^{1}L_M(\Omega) : \mbox{ tr}(u)=0 \mbox{ on } \partial\Omega\}\subset \{u\in W^{1,1}(\Omega) : \mbox{ tr}(u)=0 \mbox{ on } \partial\Omega\}=W^{1,1}_{0}(\Omega)$. So that for any $v\in\{u\in W^{1}L_M(\Omega) : \mbox{ tr}(u)=0 \mbox{ on } \partial\Omega\}$ one has $v\in W^{1,1}_{0}(\Omega)\cap W^{1}L_M(\Omega)=W^{1}_0L_M(\Omega)$.
\qed
\bibliographystyle{plain}

\end{document}